\documentclass{article}

\usepackage{fullpage}

\usepackage{graphicx}
\usepackage{graphics}
\usepackage{epsfig}
\usepackage{amsmath}
\usepackage{amssymb}
\usepackage[english]{babel}
\usepackage[latin1]{inputenc}
\usepackage{latexsym}

\usepackage{booktabs}

\def\sgn{\mbox{sign}}

\def\ds{\displaystyle}

\def\bproof{\noindent {\em Proof.} \hspace{1mm}}
\def\eproof{\hfill \ensuremath{\Box}}

\begin{document}

\title{Minimal time control of fed-batch processes\\ 
with growth functions having several maxima}

\author{A.~Rapaport${}^{1}$\footnote{A. Rapaport is with the Equipe-projet INRA-INRIA 'MERE' (Mod\'elisation Et Ressources en Eau)} and D.~Dochain${}^{2}$\footnote{Honorary Research Director FNRS, Belgium.}}
\date{\small
${}^{1}$ UMR 'MISTEA' Math\'ematiques, Informatique et STatistique pour l'Environnement et l'Agronomie (INRA/SupAgro) 2, place P. Viala, 34060 Montpellier,  France.\\ E-mail: rapaport@supagro.inra.fr\\
${}^{2}$ CESAME, Universit\'e catholique de Louvain, 4-6 avenue G. Lema\^{\i}tre 1348 Louvain-la-Neuve, Belgium.\\ E-mail: Denis.Dochain@uclouvain.be}

\maketitle

\begin{abstract}
In this paper we address the issue of minimal time optimal control of fedbatch reactor in presence of complex non monotonic kinetics, that can be typically characterized by the combination of two Haldane models. The optimal synthesis may present several singular arcs.
Global optimal trajectory results are provided on the basis of a numerical approach that considers an approximation method with smooth control inputs.
\end{abstract}

\maketitle

\section{Introduction}

Fed-batch bioreactors represent an important class of bioprocesses, mainly in the food industry (e.g. yeast production or wine making) and in the pharmaceutical industry (like the production of the vaccine against the Hepatitis B) but also e.g. for biopolymer applications (PHB). It is also very much involved in the field of enzyme production which has been developed over the past decade due to the recombinant ADN technology and via the use of filamentous micro-organisms. One of the key issues in the operation of fed-batch reactors is to optimize the process operation over a limited time period. A intensive research activity has been devoted to optimal control of (fed-batch) bioreactors mainly in the seventies and in the eighties (see e.g. \cite{O76,P85}). In this paper we address the issue of minimal time optimal control of fedbatch reactor in presence of complex non monotonic kinetics, characterized here by the combination of two Haldane models, aimed to emphasize the presence of parallel metabolic pathways to transform the limiting substrate $S$ into the biomass $B$.
For those problems, it is well knwon that the optimal synthesis is bang-bang with the possibility of a singular arc \cite{Moreno}.
For combinaisons of several non-monotonic growths,
the multiplicity of singular arcs reveals a issue for determining 
which singular arc is globally optimal.
In this work, global optimal trajectory results are provided on the basis of a numerical approach that considers an approximate problem whose optimal feedback laws are smooth.

The paper is organized as follows. In Section \ref{section-model}, we present the model and the hypotheses. In Section \ref{section-extremals}, we study the field of extremals, providing trajectories as candidate optimal solutions. Section \ref{section-approximation} presents our approximation procedure and Section \ref{section-numerical} shows how it can be applied to our problem, illustrated on two examples.
We end by a conclusion in Section \ref{section-conclusion}.

\section{The model}
\label{section-model}

We consider a reaction scheme where $n$ several bio-reactions consuming the same substrate $S$ and producing a biomass $B$ may occur simultaneously:
\begin{equation} S+B+R_{i} \longrightarrow B+B \qquad i=1,\cdots, n \end{equation}
Each reaction $i$ requires an additional resource $R_{i}$, that we shall assume to be non-limiting, and is characterized by a specific growth function $\mu_{i}(\cdot)$. Under the hypothesis of non-limitation by the auxiliary resources, one can assume that each function  $\mu_{i}(\cdot)$ depends on the concentration of the substrate $S$ only.\\

\noindent {\bf Assumption 1.} The functions $\mu_{i}(.)$ are smooth, positive away from zero and null at $0$.\\

Moreover, we shall assume that each bio-reaction (that could be associated to  different metabolic pathways) transforms the substrate into biomass with the same yield $y$. The time evolution of the concentrations $S$ and $B$ in a perfectly mixed reactor, operated in fed-batch, is given by the following dynamical system:
\begin{eqnarray}
\dot S & = & -\frac{1}{y}\sum_{i=1}^{n}\mu_{i}(S)B+\frac{Q}{V}(S_{in}-S) \label{dS} \\
\dot B & = & \sum_{i=1}^{n}\mu_{i}(S)B -\frac{Q}{V}B \label{dB}\\
\dot V & = & Q \label{dV}
\end{eqnarray} 
where $V \in (0,V_{max}]$ is the volume of the liquid phase in the tank, $S_{in}>0$ the input concentration of substrate and $Q$ the input flow rate.

In many applications (including wastewater treatment processes), a typical objective is to reach in minimal time a target:
\begin{equation}
\left\{(S,B,V)\in\mathbb{R}_{+}^{3} \, | \, S\leq S_{ref} \mbox{ and } V=V_{max}\right\}
\end{equation} 
where $S_{in}>S_{ref}>0$ via the manipulated variable $Q\in[0,Q_{max}]$.

Let us define:
\begin{eqnarray}
\mu(S) &=& \sum_{i=1}^{n}\mu_{i}(S) \\
X &=& \frac{B}{y}
\end{eqnarray}
the system dynamics (\ref{dS})-(\ref{dV}) are equivalent to the one of a single bio-reaction with specific growth rate $\mu(\cdot)$ and unitary yield factor:
\begin{eqnarray}
\label{dyn3dim}
\dot S & = & -\mu(S)X+\frac{Q}{V}(S_{in}-S \label{dS1})\\
\dot X & = & \mu(S)X -\frac{Q}{V}B \label{dX1} \\
\dot V & = & Q \label{dV1}
\end{eqnarray} 
One can easily check, from equations (\ref{dyn3dim}), that the quantity $M=V(X+S-S_{in})$ is constant along the trajectories, and consequently that it depends 
on the initial condition $(S_{0},X_{0},V_{0})$ only. By considering the number $M_{0}=V_{0}(X_{0}+S_{0}-S_{in})$, the dynamics that can be defined in the 
$(S,V)$-plane as follows:
\begin{eqnarray}
\label{dyn2dim}
\label{dynS} \dot S & = & -\mu(S)(M_{0}/V-S+S_{in})+\frac{Q}{V}(S_{in}-S)\\
\label{dynV} \dot V & = & Q
\end{eqnarray} 
with the target:
\begin{equation}
{\cal T} = [0,S_{ref}]\times\{V_{max}\}
\end{equation}
This minimal time problem has already been studied in \cite{Moreno} for cases where the function $\mu(\cdot)$ has at most one maximum. An extension with impulse control inputs has been developed in \cite{SIAM}. With the help of a clock form and Green's Theorem, Moreno  has proved that: 
\begin{enumerate}
\item  the ``bang-bang'' strategy (i.e. $Q=Q_{max}$ until $V=V_{max}$ and then $Q=0$) is optimal for any monotonic growth function $\mu(\cdot)$,
\item  the ``singular arc'' strategy (that consists in reaching and remaining at $S=\bar S$ as long as possible) is optimal 
for growth functions increasing when $S<\bar S$ and decreasing when $S>\bar S$ (and under the condition $S_{in}>\bar S>S_{ref}$).
\end{enumerate}
The clock form is used in Moreno's proof to show the global optimality of these strategies, a technique originated from the former results of Miele \cite{Miele} (see also \cite{HermesLasalle}).
For cases where the growth functions have more than one local maximum on  the interval $(0,S_{in})$, this argument can still be used but 
only for the local optimality of singular arcs.
Unfortunately, one cannot directly deduce the global optimality of the singular arc strategy.

In the present work, we shall consider the following assumption.\\

\noindent {\bf Assumption 2.} The function $\mu(\cdot)$ has a finite number of local maxima on the interval $[0,S_{in}]$. Moreover, the set
\begin{equation}
{\cal M} = \{ \bar S \in [0,S_{in}] \, \vert \, \mu(\cdot) \mbox{ is locally maximal at } \bar S \}
\end{equation}
is such that:
\begin{eqnarray}
\mbox{card}\,{\cal M} > 1\\
S_{ref}<\bar S_{-} = \min {\cal M} < \bar S_{+} = \max {\cal M} < S_{in}
\end{eqnarray}

\medskip

\begin{figure}[htbp]
\centering
\includegraphics[scale=0.5]{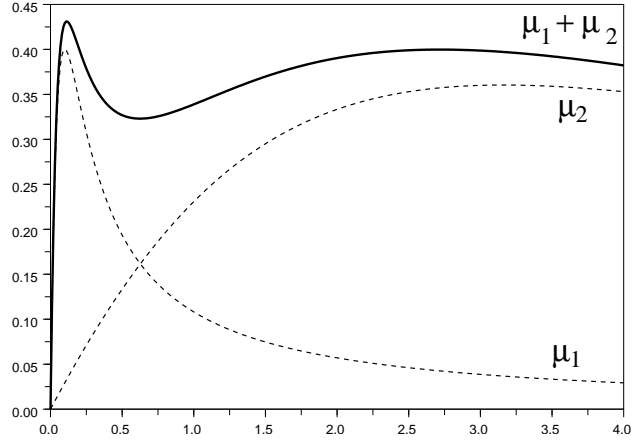}
\caption{\label{fig2haldane} Graph of the sum of two Haldane functions.}
\end{figure}

This assumption can be typically fulfilled with $n=2$ and $\mu_{1}(\cdot)$, $\mu_{2}(\cdot)$ of the Haldane type:
\begin{equation}
\mu_{i}(S)=\frac{\bar\mu_{i}S}{K_{i}+S+S^{2}/L_{i}}
\end{equation}
that admits a single maximum at $\sqrt{K_{i}L_{i}}$. 
For instance, the sum of the two Haldane functions
\[
\mu_{1}(S)=\frac{1.2\, S}{0.1+S+10\,S^{2}},\quad
\mu_{2}(S)=\frac{1.5\, S}{5+S+0.5\,S^{2}}
\]
has two local maxima (see Fig.\ref{fig2haldane}).

\section{Study of the extremals} 
\label{section-extremals}

We shall consider initial conditions on the domain ${\cal D}=[S_{ref},S_{in})\times(0,V_{max}] \ $ for our study (it sounds realistic that the initial substrate concentration is between the input 
and the desired ones).

We first assume that the maximal flow rate $Q_{max}$ is large enough for ensuring the 
local controllability of the dynamics on the domain ${\cal D}$.\\

\noindent {\bf Assumption 3. }
\begin{equation} 
Q_{max} >  \displaystyle \max_{S\in[\bar S_{-},\bar S_{+}]}\mu(S)\left(\frac{M_{0}}{S_{in}-S}+V_{max}\right)
\end{equation}

\medskip

We show now that the target ${\cal T}$ can be replaced by a punctual one.\\

\noindent {\bf Proposition 1.} \label{prop1} From any initial condition in ${\cal D}$, the optimal trajectory reaches the target ${\cal T}$ at  point $(S_{ref},V_{max})$.\\

\bproof 
Consider the curve in the $(S,V)$-plane:
\begin{equation} {\cal C}=\{ (\sigma(V),V) \, \vert \, V\in[0,V_{max}] \} \end{equation}
where $\sigma(\cdot)$ is solution of the differential equation
\begin{equation} \left\{\begin{array}{lll} \ds \frac{d\sigma}{dv} & = & \ds -\frac{\mu(S)}{Q_{max}}\left(\frac{M_{0}}{v}+S_{in}-S\right) + \frac{S_{in}-S}{v}\\ \sigma(V_{max}) & = & S_{ref} \end{array}\right. \end{equation}
Then the domain ${\cal E}$ delimited by this curve and the 
$(S,V)$ axes:
\begin{equation} 
{\cal E} = \{  (S,V) \, \vert \, 
0 \leq S \leq \sigma(V)
\mbox{ and }
\min \{ v\geq 0 \, \vert \, \sigma(v)\geq 0 \} \leq V \leq V_{max}
\; \} 
\end{equation}
is invariant whatever is the control. Assumptions 2 and 3 imply that the function $\sigma(\cdot)$  is increasing and consequently $S\leq S_{ref}$ for any $(S,V)$ in ${\cal E}$. Note also that Assumption 2 imply that the function $\mu(\cdot)$ is increasing on $[0,S_{ref}]$. From the result of Moreno, one deduces that the curve  ${\cal C}$ is an optimal trajectory. 

Consider an initial condition in ${\cal D}\setminus{\cal T}$. It does not belong to  ${\cal E}$, but if an optimal trajectory reaches the target with $S<S_{ref}$, it has to enter the domain
${\cal E}$ i.e. it has to cross the curve ${\cal C}$ with $S<S_{ref}$, which contradicts the optimality of  the curve ${\cal C}$. \eproof

\medskip

Let us write the Hamiltonian of the optimization problem:
\begin{equation} \label{hamiltonian} 
H = \lambda_{0} -\lambda_{S}\mu(S)(M_{0}/V-S+S_{in})
+ Q\left(\lambda_{S}\frac{S_{in}-S}{V} +\lambda_{V}\right) 
\end{equation}
where $\lambda_{0}\leq 0$, and the adjoint equations are:
\begin{eqnarray}
\label{dynlambdaS}
\dot \lambda_{S} & = & \lambda_{S}\left(\mu^{\prime}(S)X-\mu(S)+\frac{Q}{V}\right)\\
\dot \lambda_{V} & = & \lambda_{S}\frac{-\mu(S)M_{0}+Q(S_{in}-S)}{V^{2}}
\end{eqnarray}
with
\begin{equation} X=\frac{M_{0}}{V}+S_{in}-S \end{equation}
We define the switching function:
\begin{equation} \phi = \lambda_{S}\frac{S_{in}-S}{V}+\lambda_{V} \end{equation}
from which one deduces the optimality of ``bang-bang'' control inputs:
\begin{equation}  Q^{\star} = \left|\begin{array}{ll} 0 & \mbox{if } \phi<0\\ Q_{max} & \mbox{if } \phi>0 \end{array}\right. \end{equation}
We focus now on the characterization of singular arcs.\\

\noindent {\bf Proposition 2.} The singular arcs are trajectories
\begin{equation} \label{Vsingular}
\left\{\begin{array}{lll}
S(t) & = & \bar S\\
V(t) & = & \ds V(t_{1})e^{\mu(\bar S)(t-t_{1})}+\frac{M_{0}}{S_{in}-\bar S}\left(e^{\mu(\bar S)(t-t_{1})}-1\right)
\end{array}\right.
\end{equation}
for $t \in [t_{1},t_{2}]$,
with $t_{2}>t_{1}\geq 0$, $V(t_{1})<V_{max}$, $V(t_{2})\leq V_{max}$, where $\bar S \in {\cal M}$ and the control input is given by the feedback equation:
\begin{equation}
\label{feedbacksingular}
Q^{s}(\bar S,V)=\mu(\bar S)\left(\frac{M_{0}}{S_{in}-S}+V\right)
\end{equation}

\medskip

\bproof A singular arc strategy may occur when the switching function is identically equal to zero on a time interval of positive measure. One can easily compute:
\begin{equation} \label{dotphi} \dot \phi = \lambda_{S}\frac{S_{in}-S}{V}\mu^{\prime }(S)X \end{equation}
Note from equation (\ref{dynS}) that from any initial condition in ${\cal D}$, one has:
\begin{equation} S(t)<S_{in} \ , \quad \forall t \geq 0 \end{equation}
and consequently one has:
\begin{equation} 
X(t)=\frac{M_{0}}{V}-S(t)+S_{in} >0 \, , \quad \forall t\geq 0 \ .
\end{equation}
Note that it is not possible to reach the target from any initial condition in ${\cal D}\setminus {\cal T}$, with a constant control input $Q=0$ and $Q=Q_{max}$. 
Consequently, $\phi$ has to be equal to zero at a certain time. Note also from equation (\ref{dynlambdaS}) that the sign of $\lambda_{S}$ is constant or $\lambda_{S}$ is identically equal to $0$.
In this latter case, $\lambda_{V}$ is constant and has to be non-zero from the Maximum Principle. Yet then $\phi=\lambda_{V}$ cannot 
be equal take to zero. So $\lambda_{S}$ is never equal to zero, and when $\phi=0$ one can conclude from $H=0$, where $H$ is the Hamiltonian defined in (\ref{hamiltonian}), that $\lambda_{S}$ has to be negative.

We deduce from equation (\ref{dotphi}) that a necessary condition for  an extremal to be singular is to have:
\begin{equation} \mu^{\prime}(S)=0 \end{equation}
The zero of $\mu^{\prime}(\cdot)$ being isolated, we deduce that the singular arcs are of the form:
\begin{equation} S(t)=\bar S \mbox{ with } \mu^{\prime}(\bar S)=0 \end{equation}
Along with this condition, one can easily compute:
\begin{equation} \ddot \phi = \frac{\lambda_{S}(S_{in}-\bar S)}{V}\mu^{\prime\prime}(\bar S)X\dot S \end{equation}
Assumption 3 implies that for $\phi\neq 0$, one has:
\begin{equation} \sgn(\dot S) = \sgn(\phi) \end{equation}
Consequently, in the neighborhood of $\bar S$, one has:
\begin{equation} \sgn(\ddot \phi) = -\sgn(\mu^{\prime\prime}(\bar S))\sgn(\phi) \end{equation}
From the classification of fold points for time optimal control in the plane \cite{BonnardChyba}, one obtains that a point $(\bar S,V,\lambda_{S},\lambda_{V})$
such that $\phi=0$ is:
\begin{itemize}
\item[-] elliptic when $\mu^{\prime\prime}(\bar S)>0$, and the optimal trajectory in its neighborhood is bang-bang,
\item[-] hyperbolic when $\mu^{\prime\prime}(\bar S)<0$, and the optimal trajectory in its neighborhood can have a singular arc.
\end{itemize}
Finally, a necessary condition for an extremal to be singular is to have  $\mu^{\prime}(\bar S)=0$ and $\mu^{\prime\prime}(\bar S)<0$ which amounts to have $\bar S\in {\cal M}$. 

Having $S(t)=\bar S$ on a time interval $[t_{1},t_{2}]$ implies $\dot S=0$, and from equation (\ref{dynS}) one deduces the expression of the control given in (\ref{feedbacksingular}). 
Then from equation (\ref{dynV}) the variable $V$ is solution of the 
ordinary differential equation
\begin{equation} \dot V=\mu(\bar S)\left(\frac{M_{0}}{S_{in}-S}+V\right) \ , \end{equation} 
whose explicit solution is given in (\ref{Vsingular}).\eproof

\medskip

Proposition 2 gives only a local optimality result. Away from the set ${\cal M}\times (0,V_{max})$, we know that the optimal control is either  $0$ or $Q_{max}$ and can switch, but we do not know a priori 
\begin{itemize}
\item[-] towards which singular arc, defined by the value of $\bar S \in {\cal M}$, it is optimal to go?
\item[-] if is it optimal to quit a singular arc for reaching another singular one?
\end{itemize}
To address the global optimality, we consider now a numerical approach.

\section{An approximation procedure}
\label{section-approximation}

Solving numerically minimal time problems with dynamics that are affine w.r.t. to the control input is usually intricate when one does know a priori the switching surfaces. The shooting function based on the integration of the Hamiltonian system is usually not smooth when the optimal control is discontinuous \cite{Trelat}.  We consider here a smoothing method based on a idea originally proposed in \cite{Spinelli}.
We first define a new control:
\begin{equation} u_{1} = \frac{2Q}{Q_{max}}-1 \in [-1,1] \end{equation}
and denote:
\begin{equation} \xi = \left[\begin{array}{c} S\\ V  \end{array}\right] \end{equation}
Then one can note that the dynamics (\ref{dynS})-(\ref{dynV})  can be rewritten as follows:
\begin{equation} \label{originaldynamics} \dot \xi = F(\xi) +G_{1}(\xi)u_{1} \ , \quad \xi(0)=z_{0} \in {\cal D} \end{equation}
where
\begin{eqnarray} 
F(\xi) &=& \left[\begin{array}{c}  -\mu(S)\left(\frac{M_{0}}{V}-S+S_{in}\right)\\[2mm] 0 \end{array}\right]+G(\xi) \\
G_{1}(\xi) &=& \frac{Q_{max}}{2}\left[\begin{array}{c}  \frac{S_{in}-S}{V}\\[2mm] 1 \end{array}\right]
\end{eqnarray}
According to Proposition \ref{prop1}, we shall consider the punctual target defined by:
\begin{equation} z_{f}=\left[\begin{array}{c} S_{ref}\\ V_{max} \end{array}\right] \ . \end{equation}
The vector field $G_{1}(\cdot)$ is nowhere equal to the zero vector. Consequently there exists another vector field $G_{2}(\cdot)$ such that:
\begin{equation} \label{rankcondition} Vect\left(G_{1}(\xi),G_{2}(\xi)\right)=\mathbb{R}^{2} \ , \quad \forall \xi \in {\cal D} \end{equation}
One can also require $G_{2}(\cdot)$ to be bounded:
\begin{equation} ||G_{2}(\xi)||\leq r < +\infty\ , \quad \forall \xi \in \mathbb{R}^{2} \ . \end{equation}
Now we consider the augmented dynamics, with an additional input $u_{2}$:
\begin{equation} \label{augmenteddynamics} \dot \xi_{\epsilon} = F(\xi_{\epsilon})+G_{1}(\xi_{\epsilon})u_{1}+ \epsilon G_{2}( \xi_{\epsilon})u_{2} \ , \quad \xi(0)=z_{0} \end{equation}
with $u_{1}^{2}+u_{2}^{2}\leq 1$ and $\epsilon\neq 0$. The Hamiltonian of this new problem  is equal to:
\begin{equation} H_{\epsilon} = p_{0} + p^{t}F(\xi) + p^{t}G_{1}(\xi)u+\epsilon p^{t}G_{2}(\xi)v \end{equation}
The adjoint vector $p$ being never equal to the zero vector (by the Maximum Principle), Condition (\ref{rankcondition}) implies that:
\begin{equation} \left[\begin{array}{c} p^{t}G_{1}(\xi)\\ \epsilon p^{t}G_{2}(\xi) \end{array}\right] \neq 0 \ . \end{equation}
Consequently the Hamiltonian $H_{\epsilon}$ is uniquely maximized by the smooth control inputs:
\begin{eqnarray}
u_{1}^{\star}(\xi,p) & = & \frac{p^{t}G_{1}(\xi)}{\sqrt{[p^{t}G_{1}(\xi)]^{2}+
[\epsilon p^{t}G_{2}(\xi)]^{2}}}\\
u_{2}^{\star}(\xi,p) & = & \frac{\epsilon p^{t}G_{2}(\xi)}{\sqrt{[p^{t}G_{1}(\xi)]^{2}+
[\epsilon p^{t}G_{2}(\xi)]^{2}}}
\end{eqnarray}
We show now that the optimal trajectories for the extended dynamics converges toward an optimal trajectory of the original problem.  This convergence has been recently studied in \cite{Silva}, but the proof we propose here is different and is based on differential inclusions.\\

\noindent {\bf Proposition 3.}  Let $\epsilon_{n}$ be a monotonic sequence of numbers converging to $0$, and $\xi_{n}(\cdot)$ a sequence of optimal trajectories for $\epsilon=\epsilon_{n}$ with the same initial condition $z_{0}$ and target $z_{f}$. Then any $\bar \xi(\cdot)$ limit of a sub-sequence, also denoted $\xi_{n}$, in the following sense:
\begin{equation}
\xi_{n}(\cdot) \to \bar \xi(\cdot) \mbox{ uniformly and }
\dot\xi_{n}(\cdot) \to w(\cdot) \mbox{ weakly in } L^{1}
\end{equation}
is an optimal trajectory for the original problem. Furthermore, there exists at least one such sub-sequence.\\

\bproof Recall from Filippov's Lemma that the set of solutions of the equations (\ref{augmenteddynamics}) for measurable control inputs is exactly the set of absolutely continuous solutions of the differential inclusion:
\begin{equation} 
\dot\xi \in \Psi_{\epsilon}(\xi) = F(\xi) + \bigcup_{u\in\mathbb{B}}  \left( G_{1}(\xi)u_{1}+\epsilon G_{2}(\xi)u_{2} \right)
\end{equation} 
Note that the set-valued maps $\Psi_{\epsilon}$  are monotonic w.r.t. $\epsilon$ in the following sense:
\begin{equation} \label{psimonotonic} \epsilon<\epsilon^{\prime} \; \Longrightarrow \; \Psi_{\epsilon^{\prime}}(\xi) \subset \Psi_{\epsilon}(\xi) \quad \forall \xi \end{equation}
Let us consider an initial condition $z_{0}$ in ${\cal D}$ and a monotonic sequence of positive numbers $\epsilon_{n}$ converging to zero. As one has $\Psi_{-\epsilon}(\cdot)  = \Psi_{\epsilon}(\cdot)$, we can consider a decreasing sequence $\epsilon_{n}$
without any loss of generality.  Let us denote by $T$ and $T_{n}$ the minimal times to reach  $z_{f}$, respectively for dynamics (\ref{originaldynamics}) and (\ref{augmenteddynamics}) for $\epsilon=\epsilon_{n}$. Property (\ref{psimonotonic}) implies that   the sequence $T_{n}$ is non decreasing and bounded from above by $T$. Consequently $T_{n}$ converges to a limit, denoted $\bar T$, such that $\bar T\leq T$.

Consider now a sequence $\xi_{n}(\cdot)$ of optimal trajectories for the minimal time problem with $\epsilon=\epsilon_{n}$. These trajectories can be prolonged up to time $\bar T$ 
(taking any admissible control input on the time interval $[T_{n},\bar T]$), and are uniformly bounded on $[0,\bar T]$. According to Dunford-Pettis Theorem, there exists a sub-sequence, also denoted  $\xi_{n}(\cdot)$ such that $\dot \xi_{n}(\cdot)$ converges weakly to $v(\cdot)$ on $[0,\bar T]$. Let us then define:
\begin{equation} \bar \xi(t)=z_{0}+\int_{0}^{t} v(s)ds \ , t \in[0,\bar T] \end{equation}
By weak convergence one has $\xi_{n}(\cdot)\to \bar \xi(\cdot)$  and  $\dot \xi_{n}(\cdot)\to \dot{\bar\xi}(\cdot)$  weakly on $[0,\bar T]$. One has also:
\begin{equation}
\dot \xi_{n} \in \Psi_{0}( \xi_{n}) + \epsilon_{n}r\mathbb{B} \mbox{ a.e. } t \in [0,\bar T]
\end{equation}
By compactness of trajectories of perturbed differential inclusions (see \cite{Clarke} or \cite{Vinter}), one obtains
\begin{equation}
\dot \xi \in \Psi_{0}(\bar \xi) \mbox{ a.e. } t \in [0,\bar T]
\end{equation}
Finally, one has
\begin{equation}
\xi_{n}(T_{n})=z_{f}
\end{equation}
and from the uniform convergence and continuity of trajectories $\xi_{n}(\cdot)$, one obtains
\begin{equation}
\bar \xi(\bar T)=z_{f}
\end{equation}
Consequently $\bar \xi(\cdot)$ is an optimal trajectory for the original problem and necessarily $\bar T=T$.\eproof

\medskip

\noindent {\bf Corollary 3.}  If the original problem admits a unique optimal trajectory $\bar \xi(\cdot)$, then any sequence of optimal trajectories $\xi_{n}(\cdot)$ for $\epsilon_{n}$ a monotonic sequence of numbers converging to $0$, converges uniformly to $\bar \xi(\cdot)$, and $\dot\xi_{n}(\cdot)$ converges weakly to $\dot{\bar\xi}(\cdot)$.\\

\noindent {\em Remark.} It is shown in \cite{Silva} that the optimal control does not necessarily converge point-wise in presence of singular arcs, and may exhibit a chattering phenomenon.
Our approach here is slightly different, as we already know the  locus of singular arcs and as our aim is to address the issues raised at the end of  Section \ref{section-extremals}. This explains why we focus on the approximation of the optimal trajectories 
instead of the optimal controls.

\section{Numerical approximation}
\label{section-numerical}

In this section, we present a methodology that consists in computing the field of extremals for the regularized problem, and deducing if possible some properties of the optimal trajectories, namely answers to the questions raised at the end of Section 
\ref{section-extremals}. We illustrate this approach on two examples of growth function $\mu(\cdot)$.\\

We solve backward in time the Hamiltonian dynamics associated to the regularized problem:
\begin{eqnarray}
\dot \xi & = & F(\xi)+G_{1}(\xi)u_{1}^{*}(\xi,p)+\epsilon_{n}G_{2}(\xi)u_{2}^{*}(\xi,p) \nonumber \\
\dot p & = & -p^{t}\partial_{\xi}F(\xi)-p^{t}\partial_{\xi}G_{1}(\xi)u_{1}^{*}(\xi,p) -\epsilon_{n}p^{t}\partial_{\xi}G_{2}(\xi)u_{2}^{*}(\xi,p) \nonumber
\end{eqnarray}
from the terminal state $(\xi,p)=(z_{f},p_{f})$ with different  non zero vector $p_{f}$. Without any loss of generality,  one can choose vectors $p_{f}$ of norm equal to one i.e.:
\begin{equation}
p_{f}=\left[\begin{array}{c} \cos(\alpha) \\ \sin(\alpha)
\end{array}\right]
\end{equation}
Taking discrete values of $\alpha$ in $[0,2\pi)$, we plot the projections of each solution in the $(S,V)$-plane. Two situations may happen :

\noindent - either we fill the domain ${\cal D}$ with a set of trajectories that do not intersect. Then each trajectory is optimal for the approximated problem, and consequently is close from an optimal one of the original problem (see Example 1 below).

\noindent - either some trajectories intersect in ${\cal D}$ and only the part before the intersection (in backward time) can be optimal.
Nevertheless, this partial information might be enough to answer the questions raised at the end of Section \ref{section-extremals} (see Example 2 below).


For the vector field $G_{2}(\cdot)$, we have simply chosen a constant one (but other choices are possible):
\begin{equation} G_{2}(\xi)=\left[\begin{array}{c} 1\\ 0 \end{array}\right] \ , \quad \forall \xi \end{equation}
When $G_{2}(\cdot)$ is a constant, the adjoint equations are  independent of $\epsilon$. Then one can show, similarly to the original problem, that $p_{S}$ and $p_{V}$ are respectively negative and  positive, for extremals with $\xi(\cdot)$ in the domain ${\cal D}$. Consequently, we take values of $\alpha$ in the interval $(\pi,3\pi/2)$ only.\\

In both examples below, values of parameters of the problem are given in Table \ref{tableparameters}.\\

\begin{table}[htbp]
\caption{\label{tableparameters}Numerical simulation parameters}
\centering
\begin{tabular}{|c|c|c|c|c|c|}
\toprule
$S_{ref}$ & $V_{max}$ & $S_{in}$ & $y$ & $M_{0}$ & $Q_{max}$\\
\midrule
0.1 & 50 & 10 & 5 & 170 & 5\\
\bottomrule
\end{tabular}
\end{table}

\begin{figure}[htbp]
\begin{center}
\includegraphics[scale=0.5]{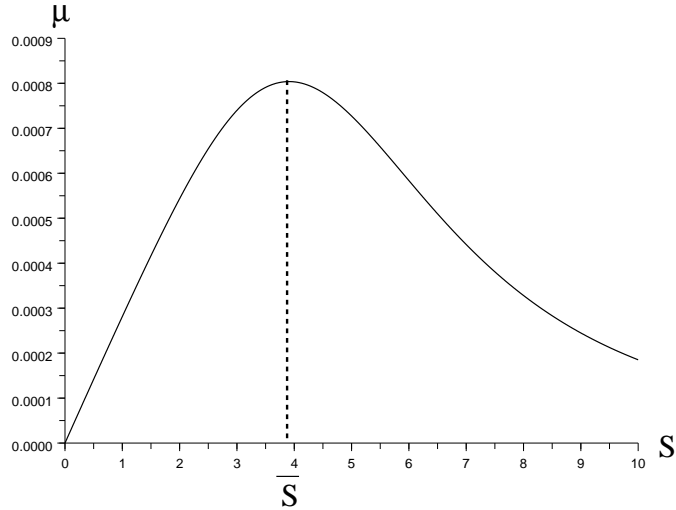}
\caption{\label{figmu1} Graph of a function $\mu(\cdot)$
with one maximum.}
\end{center}
\end{figure}

\noindent {\bf Example 1.} We first test the method for a growth function $\mu(\cdot)$ with only one maximum reached at $\bar S$ (see Fig.\ref{figmu1}), for which the optimal solution is known \cite{Moreno}.

\begin{figure}[h]
\begin{center}
\includegraphics[scale=0.5]{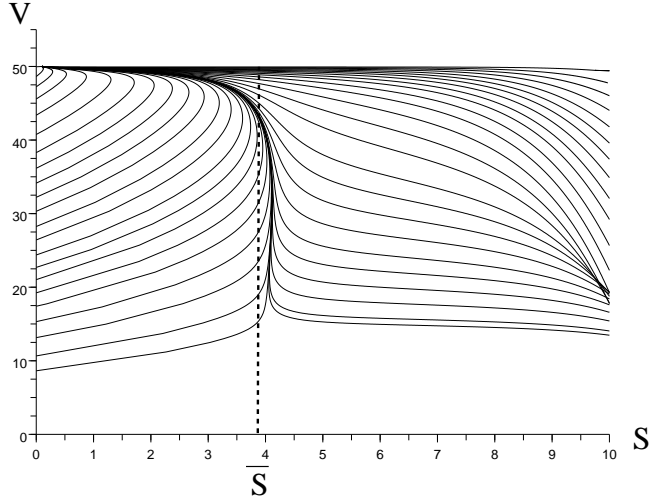}
\caption{\label{figsimu1} Extremals for the $\epsilon$-problem 
with one maximum.}
\end{center}
\end{figure}

It consists in reaching the singular arc $S=\bar S$ an staying
on this arc until reaching the boundary of the domain (i.e. $V=V_{max}$),
using the following feedback
\begin{equation}
\label{singularfeedback}
Q_{\bar S}(S,V)= \left|\begin{array}{ll}
0 & \mbox{if } S>\bar S\\
Q^{s}(\bar S,V) & \mbox{if } S=\bar S \mbox{ and } V<V_{max}\\
Q_{max}  & \mbox{if } S<\bar S \mbox{ and } V<V_{max}
\end{array}\right.
\end{equation}
Extremals for the augmented dynamics with $\epsilon=0.01$ are plotted on Fig.\ref{figsimu1}. One can see that the domain ${\cal D}$ is filled by extremals without intersection. Consequently each extremal is an optimal trajectory for the augmented dynamics, and one can check that they are close from the optimal trajectories for the original problem, given by the feedback (\ref{singularfeedback}).\\

\begin{figure}[htbp]
\begin{center}
\includegraphics[scale=0.5]{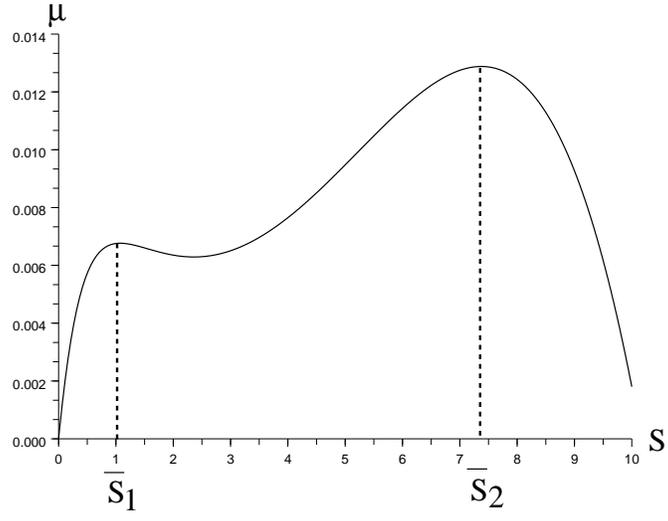}
\caption{\label{figmu2} Graph of a function $\mu(\cdot)$
with two maxima.}
\end{center}
\end{figure}

\noindent {\bf Example 2.} We consider here a growth function $\mu(\cdot)$ with two local maxima $\bar S_{1}$ and $\bar S_{2}$ (see Figure \ref{figmu2}).

Extremals for the augmented dynamics with $\epsilon=0.01$ are plotted on Fig.\ref{figsimu1}. In this case, one can see that some extremals intersect, and so one cannot conclude on their optimality on the sub-domain $B$ depicted in gray in Figure \ref{figzones}.

\begin{figure}[htbp]
\begin{center}
\includegraphics[scale=0.5]{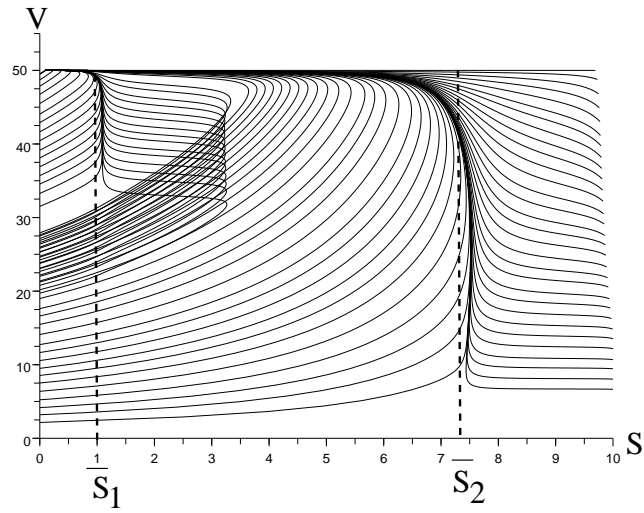}
\caption{\label{figsimu2} Extremals for the $\epsilon$-problem 
with two maxima.}
\end{center}
\end{figure}

Nevertheless, for the original problem (i.e. for $\epsilon=0$), we know that away from the singular arcs $S=\bar S_{1}$, or $S=\bar S_{2}$, the optimal control is either $0$ or $Q_{max}$ (and can switch). 
So an optimal trajectory starting from $S \in(\bar S_{1},\bar S_{2})$
 has to go toward $S=\bar S_{1}$ or $S=\bar S_{2}$ (unless it touches the boundary $V=V_{max}$) but we are not able to decide a priori 
\begin{itemize}
\item[-] if it is optimal to reach $V=V_{max}$ before reaching a singular arc,
\item[-] if not, towards which singular arc it is optimal to go,
\item[-] if it is optimal to stay on a singular arc until reaching $V=V_{max}$ (it might be better to leave one singular arc 
to go to the other one).
\end{itemize}

\begin{figure}[htbp]
\begin{center}
\includegraphics[scale=0.5]{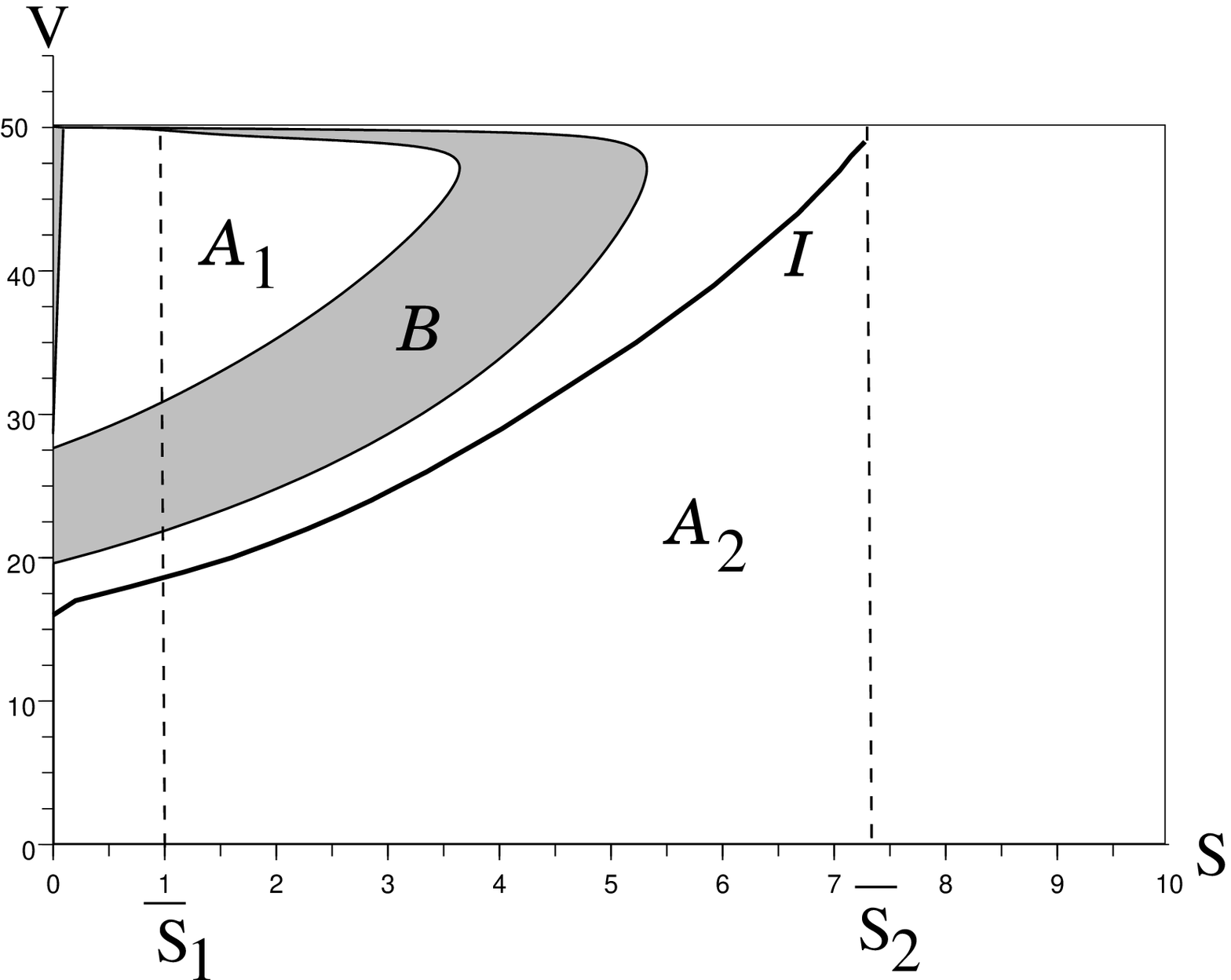}
\caption{\label{figzones} Sub-domains of optimality for reaching one of the singular arcs.}
\end{center}
\end{figure}

On Fig.\ref{figsimu1}, one can also distinguish two areas $A_{1}$ (resp. $A_{2}$), in the complementary domain of the sub-domain $B$,
depicted on Fig.\ref{figzones},
such that that extremals do not intersect and are close from the trajectories given by the feedback 
(\ref{singularfeedback}) with $\bar S=\bar S_{1}$ resp. $\bar S_{2}$.
This observation is important because it allows us to conjecture
that the optimal trajectories for the original problem reach one of the singular arcs and stay on it until reaching $V=V_{max}$.\\

\noindent {\bf Conjecture.} The optimal solution of the problem with a growth function presenting two local maxima at $\bar S_{1}$ and $\bar S_{2}$
is given by the feedback $Q_{\bar S}(\cdot)$ where
$\bar S\in\{\bar S_{1},\bar S_{2}\}$.\\

Finally, the curve $I$ on Figure \ref{figzones} has been determined numerically as the set of points for which using feedback 
$Q_{\bar S}(\cdot)$
with $\bar S=\bar S_{1}$ or  $\bar S=\bar S_{2}$ give exactly the same time for reaching the target. Then, on the left part of the domain delimited by this curve, we conjecture that 
the control $Q_{\bar S_{1}}(\cdot)$ is optimal, and $Q_{\bar S_{2}}(\cdot)$ in the right one.
Furthermore, one can observe that when $\epsilon$ get close from $0$, the boundary of the sub-domain $B$ get close from the curve $I$.

\section{Conclusion}
\label{section-conclusion}

A novel numerical method for the investigation of 
minimal time problems in the plane, that may present 
several singular arcs, has been proposed. 
It is based on an approximation 
with no singular arc and for which extremals 
can be computed straightforwardly. The method has been applied 
on the optimal control of fed-batch processes 
with non-monotonic growth function. When the field of extremals of the approximated problem has no intersection, the optimal
synthesis of the original problem can be deduced. 
Otherwise, we show that the method brings insights on optimal 
trajectories of the original problem on sub-domains of the state space.\\

\noindent {\bf Acknowledgment.} 
The authors thank Prof. C. Lobry for having pointed out the approximation procedure, and Prof. E. Tr\'elat for fruitful discussions. 
This paper presents research
results of the Belgian Programme on Inter-University Poles of
Attraction initiated by the Belgian State, Prime Minister's
office for Science, Technology and Culture. The scientific
responsibility rests with its authors.

\end{document}